\documentstyle[titlepage,twoside,12pt]{article}
\begin{document}

\pagenumbering{arabic}
\setcounter{page}{1}

\pagenumbering{arabic}

{\Large \bf Can there be a general nonlinear PDE} \\
{\Large \bf theory for existence of solutions ?} \\ \\

{\it Elem\'{e}r E Rosinger \\ Department of Mathematics \\ University of Pretoria \\ Pretoria, 0002 South Africa \\
e-mail : eerosinger@hotmail.com} \\ \\

{\bf Abstract} \\

Contrary to widespread perception, there is ever since 1994 a unified, general type
independent theory for the existence of solutions for very large classes of nonlinear systems
of PDEs. This solution method is based on the Dedekind order completion of suitable spaces of
piece-wise smooth functions on the Euclidean domains of definition of the respective PDEs. The
method can also deal with associated initial and/or boundary value problems. The solutions
obtained can be assimilated with usual measurable functions or even with Hausdorff continuous
functions on the respective Euclidean domains. \\
It is important to note that the use of the order completion method does {\it not} require any
monotonicity condition on the nonlinear systems of PDEs involved. \\
One of the major advantages of the order completion method is that it {\it eliminates} the
algebra based dichotomy "linear versus nonlinear" PDEs, treating both cases with equal ease.
Furthermore, the order completion method does {\it not} introduce the dichotomy "monotonous
versus non-monotonous" PDEs. \\
None of the known functional analytic methods can exhibit such a performance, since in addition
to topology, such methods are significantly based on algebra. \\ \\

{\bf 1. A Sample of customary perception} \\

The 2004 edition of the Springer Universitext book "Lectures on PDEs" by V I Arnold, starts
on page 1 with the statement :

\begin{quote}
"In contrast to ordinary differential equations, there is {\it no unified theory} of partial
differential equations. Some equations have their own theories, while others have no theory at
all. The reason for this complexity is a more complicated geometry ..." (italics added)
\end{quote}

The 1998 edition of the book "Partial Differential Equations" by L C Evans, starts his
Examples on page 3 with the statement :

\begin{quote}
"There is no general theory known concerning the solvability of all partial differential
equations. Such a theory is {\it extremely unlikely} to exist, given the rich variety of
physical, geometric, and probabilistic phenomena which can be modeled by PDE. Instead,
research focuses on various particular partial differential equations ..." (italics added)
\end{quote}

And yet, in 1994, in Oberguggenberger \& Rosinger, MR 95k:35002, precisely such a general
theory of existence of solutions for very large classes of nonlinear PDEs was published. For
latest developments, see Anguelov \& Rosinger. \\

In the sequel, we present the main ideas and motivations which underlie the order completion
method. The detailed mathematical developments can be found in the references, among them, in
Oberguggenberger \& Rosinger and Anguelov \& Rosinger. \\

It is on occasion worth recalling that we all do mathematics based on certain underlying ideas
and motivations. What happens is that we may hold to them for longer, and do so long enough,
so that many of them may become rather automatic. And once that happens, we do no longer - and
in fact, can no longer - review them, and do so at least now and then. This is, then, how
perceptions are established, and we end up being subjected to them. \\
Here an attempt is made to go beyond such perceptions in the realms of solving PDEs. And since
perceptions are inevitably formulated in a "meta-language" - in this case "meta" with respect
to the usual formal mathematical texts - much of what follows has to go along with that. \\ \\

{\bf 2. The class of nonlinear systems of PDEs solved} \\

In Oberguggenberger \& Rosinger it was show how to obtain solutions $U$ for all systems
of nonlinear PDEs with associated inital and/or boundary value problems, where the equations
are of the form \\

(2.1) $~~~ F(x, U(x), ~.~.~.~ , D^p_x U(x), ~.~.~.~ ) ~=~ f(x),~
                    x \in \Omega \subseteq {\bf R}^n,~  | p | \leq m $ \\

Here $F$ is any function jointly continuous in all its arguments, the right hand term $f$ can
belong to a class of discontinuous functions, the order $m \in {\bf N}$ is given arbitrary,
while the domain $\Omega$ can be any bounded or unbounded open set in ${\bf R}^n$. \\
In fact, even the functions $F$ defining the nonlinear partial differential operators in the
left hand terms of (2.1) can have certain types of discontinuities. \\

The solutions $U$ obtained in Oberguggenberger \& Rosinger can be assimilated with usual
measurable functions on the respective domains $\Omega$. \\
Recently, in Anguelov \& Rosinger this general regularity result was further improved as it
was shown that the solutions $U$ can be associated with Hausdorff continuous functions on the
same domains $\Omega$. \\

Here it is important to note that {\it Hausdorff continuous} functions are not much unlike
usual real valued continuous functions. Indeed, on suitable {\it dense} subsets of their
domains of definition, Hausdorff continuous functions have as values real numbers, and are
completely {\it determined} by such values. On the rest of their domains of definition,
Hausdorff continuous functions can have values given by bounded or unbounded closed intervals
of real numbers. Also, every real valued function which is continuous in the usual sense will
be Hausdorff continuous as well. \\

One of the major advantages of the order completion method is that it {\it eliminates} the
dichotomy between linear, and on the other hand, nonlinear PDEs, treating both cases with
equal ease. Indeed, the dichotomy between linear and nonlinear follows from the vector space
structure of the spaces of functions on which the partial differential operators act. In this
way, this dichotomy is of an algebraic nature. On the other hand, partial orders are more
basic mathematical structures, and as such, they do not, and simply cannot, differentiate
between linear and nonlinear. \\
Clearly, functional analytic methods, which rely not only on topological but also algebraic
structures cannot exhibit such a performance, since are bound to discriminate between linear
and nonlinear, see details in section 9. \\ \\

{\bf 3. A short history of difficulties in solving linear and} \\
\hspace*{0.5cm} {\bf nonlinear PDEs} \\

The first general, that is, {\it type independent} existence result for solutions of rather
arbitrary nonlinear systems of PDEs was obtained in 1874, when upon the suggestion of K
Weierstrass, Sophia Kovalevskaia gave a rigorous proof for an earlier theorem of Cauchy,
published in the 1821, in his Course d'Analyse. This result although completely general as far as the type independent
nonlinearities involved ar concerned, assumes however, that in the systems of PDEs of the form
(2.1) both $F$ and $f$ are analytic. In addition one also assumes initial value problems on
non-characteristic analytic hypersurfaces, while boundary value problems are not treated by
the respective Cauchy-Kovalevskaia theorem. \\
However, in such a highly particular situation concerning the regularity of the PDEs and the
data involved, the solutions obtained are proved to exist always, and also to be unique and
analytic. \\
The problem with that classical {\it existence, uniqueness} and {\it regularity} result is
that, typically for nonlinear PDEs, such analytic solutions do not - and in general, cannot -
exist globally on the whole of the domain of the respective PDEs, but only in certain
neghbourhoods of the analytic hypersurfaces on which the initial values are given. This is,
therefore, not due to the specific method of proof of Kovalevskaia. \\

Here however it is important to note that the failure of the existence of global analytic
solutions is but a part of a far more general phenomenon, since even linear, let alone
nonlinear PDEs may fail to have smooth, or even merely classical solutions, even in the case
of solutions of major applicative interest. After all, even linear constant coefficient PDEs
have nonclassical solutions of particular interest, such as those given by Green functions. In
the nonlinear case, difficulties start with the simple ODE given by $U_t = U^2$  which does
not have global classical solutions either, except for the trivial solution $U = 0$. As for
the nonlinear shock wave equation $U_t + U U_x = 0$, its nonclassical solutions are precisely
those which model shocks. \\

What may be interesting, and also worthwhile to note with respect to the mentioned
Cauchy-Kovalevskaia theorem, are the following three facts :

\begin{itemize}

\item The rigorous proof by Kovalevskaia of that theorem on solutions of general nonlinear
systems of PDEs predates by about two decades the corresponding general theorem on solving
systems of nonlinear ODEs defined by continuous expressions. Indeed, the existence of
solutions for such ODEs was given by Charles Emile Picard in his 1894 Comptes Rendu Acad. Sci.
Paris paper, where the associated Cauchy problem was solved by the method of successive
approximations.

\item The only so called "hard" mathematics used in the proof of the Cauchy-Kovalevskaia
theorem is the formula for the summation of a convergent geometric progression, the rest of
the proof being but a succession of rather elementary, even if quite involved, estimates of
terms in power series. In this way, the proof of the Cauchy-Kovalevskaia theorem does {\it
not} involve methods of functional analysis, and certainly it could {\it not} involve such
methods at the time in the 1880s when it was given. On the other hand, the proof of the
corresponding general existence result for solutions of nonlinear systems of ODEs does involve
a fixed point argument in suitable spaces of functions which are complete in their respective
topologies.

\item The result in the Cauchy-Kovalevskaia theorem - when considered on its own original
terms of type independent nonlinear generality - could {\it not} so far be improved in those
very terms, regardless of all the advances in functional analysis of the last more than a
century. The only such improvement of the classical result in the Cauchy-Kovalevskaia theorem
was obtained in 1985, without however using functional analytic methods, see section 4. Indeed,
when it comes to {\it type independent nonlinear generality}, the functional analytic methods
used in solving PDEs could only bring about improvements - and often quite dramatic ones - only
in a variety of far more particular cases than the type independent nonlinear generality dealt
with in the Cauchy-Kovalevskaia theorem. In this way, in spite of more than one century of
functional analysis, the classical Cauchy-Kovalevskaia theorem still remains a {\it maximal}
result, except for its extension mentioned in section 4, which does not use functional
analysis.

\end{itemize}

In the early 1950s, soon after the introduction of the linear theory of distributions by L
Schwartz, it was proved independently by Malgrange and Ehrepreis that in case of a single PDE
of the form (2.1), if the left hand term $F$ is linear and with constant coefficients, while
$f$ is the Dirac delta distribution, then (2.1) always has a global, so called, fundamental
solution given by a suitable Schwartz distribution. \\

This rather general linear result appeared to suggest that a similar result could be obtained
in the more general case when $F$ in (2.1) is linear and with smooth coefficients. L Schwartz
himself is known to have conjectured such a generalization, and furthermore, as it appears, he
suggested it at the time to Francois Treves as a subject for his doctoral thesis. \\

However, in 1957, Hans Lewy showed that the rather simple linear first order PDE in three
space variables and with first degree polynomial coefficients \\

(3.1) $~~~ (D_x + i D_y - 2(x + y)D_z) U(x, y, z) ~=~ f(x, y, z), \\ \\
                           \hspace*{10cm}  (x, y, z) \in {\bf R}^3 $ \\

does not have any Schwartz distribution solutions in any neighbourhood of any point in
${\bf R}^3$, for a large class of smooth right terms $f$. In 1967, Shapiro gave a similar
example of a smooth linear PDE which does not have solutions in Sato's hyperfunctions. \\

In the early 1960s, L H\"{o}rmander gave certain {\it necessary} conditions for the
solvability in distributions of arbitrary linear smooth coefficient PDEs, see Rosinger
[7, pp. 37-39], [8, pp. 212-214]. \\

Regarding the perception that nontrivial general, type independent results are just about
impossible to obtain related to PDEs, it is worth noting that the Malgrange-Ehrenpreis result
on fundamental solutions is precisely such a nontrivial general and type independent existence
result within the range of all linear and constant coefficient PDEs. \\
The necessary condition for the existence of distributions solutions given by H\"{o}rmander is
also a nontrivial general and type independent result, this time within the much larger class
of all linear smooth coefficient PDEs. \\ \\

{\bf 4. Nonlinear algebraic theory of generalized functions} \\

This nonlinear theory - see 46F30 in the AMS Subject Classification 2000 at
www.ams.org/index/msc/46Fxx.html - was started in the 1960s, Rosinger [1-17], and is based on
the construction of all possible differential algebras of generalized functions which contain
the Schwartz distributions. That theory has managed to come quite near to solving the Lewy
impossibility. Yet it did not solve it completely, although it obtained generalized function
solutions for large classes of linear and nonlinear PDEs. As an example, back in 1985, it
obtained the first {\it global} existence result for the general nonlinear PDEs in the
classical Cauchy-Kovalevskaia theorem. And the respective global solutions are analytic on the
whole of the domain of the PDEs, except for certain closed and nowhere dense subsets, which can
be chosen to have zero Lebesgue measure, see Rosinger [7, pp. 259-266], [8, pp. 101-122],
[9]. \\ \\

{\bf 5. The order completion method} \\

Surprisingly, the order completion method in solving general nonlinear systems of PDEs of the
form (2.1) is based on certain very simple, even if less than usual, approximation properties,
see Oberguggenberger \& Rosinger [pp. 12-20]. To give here an idea about the ways the order
completion method works, we mention some of these approximations here in the case of one
single nonlinear PDE of the form (2.1). \\

Let us denote by $T(x, D)$ the left term in (2.1), then we have the basic approximation
property : \\

{\bf Lemma 5.1} \\

$ \begin{array}{l}
                 \forall~~ x_0 \in \Omega,~~ \epsilon > 0 ~~: \\ \\
                 \exists~~ \delta > 0,~~ P ~~\mbox{polynomial in}~~ x \in
                                                                {\bf R}^n~~ : \\ \\
                 ~~~ | | x - x_0 | | ~\leq~ \delta ~~\Longrightarrow~~ f(x) - \epsilon ~\leq~
                                                T(x,D) P(x) ~\leq~ f(x)
   \end{array}$ \\

\hfill $\Box$ \\

Consequently, we obtain : \\

{\bf Proposition 5.1} \\

$~~~ \begin{array}{l}
              \forall~~ \epsilon > 0~~ : \\ \\
              \exists~~ \Gamma_\epsilon \subset \Omega ~~\mbox{closed,~ nowhere~
                  dense~ in}~~  \Omega,~~ U_\epsilon \in C^\infty (\Omega)~~ : \\ \\
              ~~~ f - \epsilon ~\leq~ T(x, D) P ~\leq~ f ~~\mbox{on}~~ \Omega \setminus
                                                             \Gamma_\epsilon
   \end{array} $ \\

Furthermore, one can also assume that the Lebesgue measure of $\Gamma_\epsilon$ is zero,
namely \\

$~~~ mes~ (\Gamma_\epsilon) ~=~ 0$. \\

\hfill $\Box$ \\

Let us now note that, see Anguelov, or Anguelov \& Rosinger \\

$~~~ C^0 (\Omega) \subset {\bf H} (\Omega) $ \\

and the set ${\bf H} (\Omega)$ of {\it Hausdorff continuous} functions on $\Omega$ is
Dedekind order complete. \\
Consequently, we obtain the following basic result on the {\it existence} and {\it regularity}
of solutions for nonlinear PDEs of the form (2.1) : \\

{\bf Theorem 5.1} \\

$~~~ T(x, D) U(x) ~=~ f(x),~~ x \in \Omega $ \\

has solutions $U$ which can be assimilated with Hausdorff continuous functions, for a class of
discontinuous functions $f$ on $\Omega$, class which contains the continuous functions on
$\Omega$. \\

\hfill $\Box$ \\

We give here some more details related to Theorem 5.1 above. In view of Proposition 5.1, we
shall be interested in spaces of piecewise smooth functions given by \\

(5.1) $~~~ C^l_{nd} (\Omega) ~=~
      \left \{ ~~ u ~~
          \begin{array}{|l}
              \exists ~~ \Gamma \subset \Omega~~ \mbox{closed, nowhere dense}~~ : \\ \\
               ~~~~~ *)~ u : \Omega \setminus \Gamma \rightarrow {\bf R} \\ \\
               ~~~~ **)~ u \in C^l (\Omega \setminus \Gamma)
               \end{array} ~\right \} $ \\ \\

where $l \in {\bf N}$. It is easy to see that we have the inclusions \\

(5.2) $~~~ T(x, D)~ C^m_{nd} (\Omega)~\subseteq~ C^0_{nd} (\Omega) ~\subset~
                                                                     {\bf H} (\Omega)$ \\

In this way, we obtain the following more precise formulation of the result in Theorem 5.1 on
the existence and regularity of solutions : \\

{\bf Theorem 5.1*} \\

(5.3) $~~~ T(x, D)^{\#}~ ( C^m_{nd} (\Omega) )^{\#}_T ~=~ ( C^0_{nd} (\Omega) )^{\#}
                                                            ~\subset~ {\bf H} (\Omega) $ \\

\hfill $\Box$ \\

Here $( C^m_{nd} (\Omega) )^{\#}_T$ and $( C^0_{nd} (\Omega) )^{\#}$ are Dedekind order
completions of $C^m_{nd} (\Omega)$ and $C^0_{nd} (\Omega)$, respectively, when these latter
two spaces are considered with suitable partial orders. The respective partial order on
$C^m_{nd} (\Omega)$ may depend on the nonlinear partial differential operator $T(x, D)$ in
(5.2), while the partial order on $C^0_{nd} (\Omega)$ is the natural point-wise one at the
points where two functions compared are both continuous. \\
The operator $T(x, D)^{\#}$ is a natural extension of the nonlinear partial differential
operator $T(x, D)$ in (5.2) to the mentioned Dedekind order completions. \\

The meaning of (5.3) is twofold : \\

\begin{itemize}

\item for every right hand term $f \in ( C^0_{nd} (\Omega) )^{\#}$ in (2.1), there exists a
solution $U \in ( C^m_{nd} (\Omega) )^{\#}_T$, and the set $( C^0_{nd} (\Omega) )^{\#}$
contains many discontinuous functions beyond those piecewise discontinuous ones, see
Oberguggenberger \& Rosinger,

\item the solutions $U$ can be assimilated with Hausdorff continuous functions on $\Omega$,
see Anguelov \& Rosinger.

\end{itemize}

\bigskip

{\bf 6. Comparison with methods in Functional Analysis} \\

The order completion method is a powerful {\it alternative} to the usual functional analytic
ones, when solving linear or nonlinear PDEs. Details in this regard are presented in
Oberguggenberger \& Rosinger [chap. 12]. Certainly, the order completion method is {\it not}
meant to replace the functional analytic ones, the latter being useful in obtaining stronger
results in a large variety of particular PDEs. \\
Here, we shall only mention the following. Functional analytic methods in solving PDEs are
based on the {\it topological completion} of uniform spaces, such a normed or locally convex
vector spaces of suitably chosen functions. In this respect, the comparative advantages of the
order completion method can shortly be formulated as follows :

\begin{itemize}

\item unlike the functional analytic methods, which are geared more naturally to the solution
of linear PDEs, the order completion method performs equally well in the case of both linear
and nonlinear PDEs, see section 9 below,

\item unlike the functional analytic methods, which face considerable difficulties when
dealing with initial, and especially, boundary value problems, the order completion method
performs without significant additional troubles in such situations,

\item the order completion method gives solutions which can be assimilated with usual
measurable, or even Hausdorff continuous functions, and thus the solutions obtained are not
merely distributions, generalized functions or hyperfuntions.

\end{itemize}

As an illustration of the comparative situation in these two methods let us consider on a
bounded Euclidean domain $\Omega$, which has a smooth boundary $\partial \Omega$, the
following well known linear boundary value problem \\

(6.1) $~~~ \begin{array}{l}
                    \Delta~ U(x) ~=~ f(x),~~~ x \in \Omega \\ \\
                     U ~=~ 0 ~~~\mbox{on}~~ \partial \Omega
             \end{array} $ \\

As is well known, for every given $f \in C^\infty ( \bar \Omega )$, this problem has a unique
solution $U$ in the space \\

(6.2) $~~~ X ~=~ \left \{~~ v \in C^\infty (\bar \Omega) ~~|~~ v
~=~ 0 ~~\mbox{on}~~
                                                            \partial \Omega ~~\right \} $ \\

It follows that the mapping \\

(6.3) $~~~ X \ni v ~\longmapsto~ | | \Delta v | |_{L^2 (\Omega)} $ \\

\bigskip
defines a norm on the vector space $X$. Now let \\

(6.4) $~~~ Y ~=~ C^\infty (\bar \Omega) $ \\

be endowed with the topology induced by $L^2 (\Omega)$. Then in view of (6.1) - (6.4), it
follows that the mapping  \\

(6.5) $~~~ \Delta : X ~\rightarrow~ Y $ \\

is a uniform continuous linear bijection. Therefore, it can be extended in a unique manner to
an isomorphism of Banach spaces \\

(6.6) $~~~ \Delta : \bar X ~\rightarrow~ \bar Y ~=~  L^2 (\Omega) $ \\

In this way one has the classical existence and uniqueness result \\

(6.7) $~~~ \begin{array}{l}
                  \forall~~~ f \in L^2 (\Omega)~~ : \\ \\
                  \exists~!~~ U \in \bar X ~~: \\ \\
                   ~~~~~~ \Delta U = f
              \end{array} $ \\

The power and simplicity - based on linearity and topological completion of uniform spaces -
of the above classical existence and uniqueness result is obvious. This power is illustrated
by the fact that the set $\bar Y =  L^2 (\Omega)$ in which the right hand terms $f$ in (6.1)
can now be chosen is much {\it larger} than the original $Y = C^\infty (\bar \Omega)$.
Furthermore, the existence and uniqueness result in (6.7) does not need the a priori knowledge
of the structure of the elements $U \in \bar X$, that is, of the respective generalized
solutions. This structure which gives the regularity properties of such solutions can be
obtained by a further detailed study of the respective differential operators defining the
PDEs under consideration, in this case, the Laplacian $\Delta$. And in the above specific
instance we obtain \\

(6.8) $~~~ \bar X ~=~ H^2 (\Omega) \cap H^1_0 (\Omega) $ \\

As seen above, typically for the functional analytic methods, the generalized solutions are
obtained in topological completions of vector spaces of usual functions. And such completions,
like for instance the various Sobolev spaces, are defined by certain linear partial
differential operators which may happen to {\it depend} on the PDEs under consideration. \\

In the above example, for instance, the topology on the space $X$ obviously {\it depends} on
the specific PDE in (6.1). Thus the topological completion $\bar X$ in which the generalized
solutions $U$ are found according to (6.7), does again {\it depend} on the respective PDE. \\

On the other hand, with the method of order completion we are {\it no longer} looking for
generalized solutions, and instead, a {\it type independent} and {\it universal} or {\it
blanket} regularity property is attained, since the solutions obtained can always be
assimilated with usual measurable functions, or even with Hausdorff continuous functions.
Similar to the functional analytic methods, however, the order completion method obtains the
solutions in spaces which may again be related to the specific nonlinear partial differential
operators $T(x, D)$ in the equations of form (2.1). \\ \\

{\bf 7. Solving equations by extending their domains} \\
\hspace*{0.5cm} {\bf of definition : the three classical methods} \\

The ancient case of {\it solving an equation}, which shocked Pythagoras two and a half
millennia ago, is given by \\

(7.1) $~~~ x^2 ~=~ 2 $ \\

This is of the general form \\

(7.2) $~~~ E(x) ~=~ c $ \\

where we are given a mapping \\

(7.3) $~~~ E : X \rightarrow Y $ \\

together with a specific $c \in Y$, and then we want to find a solution $x \in X$ so that
(7.2) holds. \\

What shocked Pythagoras was that (7.1) could not be solved if one restricted oneself to $X =
{\bf Q}$ in (7.3). And it took no less than about two millennia or more, till we could
rigorously extend $X ={\bf Q}$ to $\bar X = {\bf R}$, and thus obtain a well defined solution
$x = \pm \sqrt 2$ of (7.1). \\

In this way, ever since, we have the following model lesson in solving equations :

\begin{itemize}

\item if one cannot solve (7.2) within the framework of (7.3), then one can try to solve it in
the {\it extended} framework \\

\end{itemize}

(7.4) $~~~ \bar E : \bar X \rightarrow \bar Y $ \\

where $X \subset \bar X$, $Y \subseteq \bar Y$, and $\bar E$ is such that we have the
commutative diagram \\

\begin{math}
\setlength{\unitlength}{0.2cm}
\thicklines
\begin{picture}(60,21)

\put(11,19){$X$}
\put(28,21){$E$}
\put(16,19.5){\vector(1,0){28.5}}
\put(47,19){$Y$}
\put(0,11){$(7.5)$}
\put(12,17){\vector(0,-1){11.5}}
\put(13.5,11){$\subseteq$}
\put(48,17){\vector(0,-1){11.5}}
\put(44.7,11){$\subseteq$}
\put(10,2){$~~ \bar X$}
\put(16,2.6){\vector(1,0){28.5}}
\put(47,2){$~ \bar Y$}
\put(28,-0.5){$\bar E$}

\end{picture}
\end{math} \\

Here however, we face the following problems :

\begin{itemize}

\item how to choose or construct $\bar X$, and then how to interpret the new, or so called
{\it generalized} solutions $x \in \bar X \setminus X$, which two questions altogether
constitute but the celebrated regularity problem,

\item how to do the same for $\bar Y$, which nevertheless need not always be done, since we
can often stay with $Y$ in (7.4) and only have to extend $X$ to $\bar X$,

\item how to define the extension $\bar E$, which often, and typically in the nonlinear case,
is not a trivial problem.

\end{itemize}

Fur further detail, we can now recall that with the equation (7.1) we had to \\

(7.6) $~~~ go~ from~~~ {\bf Q}~~~ to ~~~{\bf \bar Q} ~=~ {\bf R} $ \\

On the other hand, with the equation \\

(7.7) $~~~ x^2 + 1 ~=~ 0 $ \\

we had to \\

(7.8) $~~~ go~ from~~~ {\bf R}~~~ to ~~~{\bf C} $ \\

However, there is a vast difference between (7.6) and (7.8), respectively, between solving
(7.1) and (7.7). Indeed, we solve (7.7) through the extension (7.8) which is a mere {\it
algebraic adjoining} of an element, in this case, of $i = \sqrt -1$ to ${\bf R}$. \\

On the other hand, when solving (7.1), the extension (7.6) can be seen as at least {\it three
different}, even if in this particular case equivalent, constructions, namely, through :

\begin{itemize}

\item topology

\item  algebra

\item order.

\end{itemize}

And to be more precise, we have :

\begin{itemize}

\item The Cauchy-Bolzano method is {\it ring theoretic} plus {\it topological}, and it is
applied to ${\bf Q}$, as it obtains ${\bf R}$ according to the quotient construction in
algebras \\

(7.9) $~~~{\bf R} ~=~ {\cal A} / {\cal I} $ \\

where ${\cal A} \subset {\bf Q}^{\bf N}$ is the algebra of Cauchy sequences of rational
numbers, while ${\cal I}$ is its ideal of sequences convergent to zero.

\item The method of Dedekind is based on the {\it order completion} of ${\bf Q}$.

\end{itemize}

However, the Cauchy-Bolzano method can be generalized in two directions :

\begin{itemize}

\item In the topological generalization the algebraic part can be omitted, and instead, one
only uses the {\it topological completion} of uniform spaces, here of the usual metric space
on ${\bf Q}$.

\item In the algebraic generalization it is possible to extract the abstract essence of (7.9),
and simply start with a suitable algebra ${\cal A}$, and an ideal ${\cal I}$ in it. Such a
construction can indeed be rather abstract, since it need {\it not} involve any topology on
${\bf Q}$ or ${\bf R}$, as it happens for instance, when constructing the nonstandard reals
$^*{\bf R}$, namely \\

$~~~ ^*{\bf R} ~=~ {\cal A} / {\cal I} $ \\

Here one takes ${\cal A} = {\bf R}^{\bf N}$, that is, the algebra of all sequences of real
numbers, while the ideal ${\cal I}$ is defined by any given free ultrafilter on ${\bf N}$. \\
A rather general version of such an abstract approach, which however makes a certain limited
use of topology, has been introduced and extensively used in the nonlinear algebraic theory of
generalized functions under the AMS classification index 46F30, as mentioned in section 4
above.

\end{itemize}

What is done in the method in 46F30 is to generalize the Cauchy-Bolzano method by retaining its
ring theoretic algebraic aspect, while the topological one is weakened to the certain extent of
being confined to the topologies of Euclidean spaces only. \\

What is done in the method introduced in Oberguggenberger \& Rosinger, and further developed in
Anguelov \& Rosinger, is the extension of the classical Dedekind order completion method, used
in the construction of ${\bf R}$ from ${\bf Q}$, to suitable spaces of piece-wise smooth
functions. \\

An important fact to note is that both the topological and order completion methods give us the
property that

\begin{itemize}

\item  ${\bf Q}$ is {\it dense} in ${\bf R}$

\end{itemize}

in the respective sense of topology or order. In this way, the elements of the extension of
${\bf Q}$, that is, the elements of ${\bf R}$, are in the corresponding sense {\it arbitrarily
near} to the elements of the extended space ${\bf Q}$. Thus the elements in the extension can
arbitrarily be {\it approximated} by elements of the extended space, be it in the sense of
topology, or respectively, order. \\
Furthermore, both through the methods of topology and order, one obtains ${\bf R}$ in a {\it
unique} manner, up to a respective isomorphism. \\

In this way both the topological and order completion methods have the {\it double} advantage
that

\begin{itemize}

\item the elements of the extension are not too strange conceptually,

\end{itemize}

and furthermore

\begin{itemize}

\item the elements of the extension are near to elements of the extended space, within
arbitrarily small error.

\end{itemize}

This {\it density} property remains also in the general Dedekind order completion method used
in Oberguggenberger \& Rosinger and Anguelov \& Rosinger. Indeed, in (5.3) we have that
$C^m_{nd} ( \Omega )$ and $C^0_{nd} ( \Omega )$ are order dense in $( C^m_{nd} ( \Omega ) )^{\#}_T$ and
$( C^0_{nd} ( \Omega ) )^{\#}$, respectively. \\

Connected with the general extension method in (7.3) - (7.5) one can note that, on occasion,
the following {\it convenient} situation may occur : the extended mapping \\

(7.10) $~~~ \bar E : \bar X ~~\longrightarrow~~ \bar Y $ \\

may turn out to be an {\it isomorphism} of the respective algebraic, topological or order
structures used on $X$ and $Y$, when constructing the corresponding extensions $\bar X$ and
$\bar Y$. In such a case, and when one has a better understanding of the structure of the
elements in $\bar Y$, one can obtain in addition a {\it regularity} type result concerning the
so called generalized solutions $x \in \bar X \setminus X$ of the equations (7.2), since such
generalized solutions can be {\it assimilated} - through the isomorphism $\bar E$ - with the
corresponding elements $\bar E ( x ) \in \bar Y$. \\

A classical example of such an isomorphism (7.10) happens, for instance, in (6.5), (6.6), when
the boundary value problem (6.1) is solved by using well known functional analytic methods. \\
In that specific instance, however, the suitable further use of functional analytic methods
can lead to the {\it additional} regularity property of generalized solutions in $\bar X$, as
given in (6.8). Nevertheless, the Banach space isomorphism (6.6) - which in that case is but
the particular form taken by (7.10) - is in itself already a {\it first} regularity result
about the structure of the elements of $\bar X$. \\

The above convenient situation of an isomorphism of type (7.10) can appear as well when using
the order completion method in solving nonlinear system of PDEs. This is the reason why the
solutions obtained in Oberguggenberger \& Rosinger could be assimilated with usual measurable
functions, while in Anguelov \& Rosinger, they can be assimilated with the much {\it more
regular} Hausdorff continuous functions. \\
More specifically, in (5.3), the extended mappings $T ( x, D )^{\#}$ prove to be {\it order
isomorphisms} between the spaces $C^m_{nd} ( \Omega )^{\#}_T$ and
$C^0_{nd} ( \Omega )^{\#}$. \\
This is then, in essence, the reason why the solutions of nonlinear systems of PDEs of the
form (2.1) could earlier be assimilated with usual measurable functions, and can now be
assimilated with Hausdorff continuous functions. \\ \\

{\bf 8. The need for extensions in the case of solving PDEs}. \\

Let us now associate with each nonlinear PDE in (2.1) the corresponding nonlinear partial
differential operator defined by the left hand side, namely \\

(8.1)~~~ $ T ( x, D ) U ( x ) ~=~ F ( x, U ( x ), ~.~.~.~ , D^p_x U ( x ), ~.~.~.~ ),~~~
                                        x \in \Omega $ \\

{\it Two} facts about the nonlinear PDEs in (2.1) and the corresponding nonlinear partial
differential operators $T ( x, D )$ in (8.1) are important and immediate

\begin{itemize}

\item The operators $T ( x, D )$ can {\it naturally} be seen as acting in the {\it classical}
context, namely

\end{itemize}

(8.2)~~~ $ T ( x, D ) ~:~ {\cal C}^m ( \Omega ) \ni U ~~\longmapsto~~ T ( x, D ) U \in
                                          {\cal C}^0 ( \Omega ) $ \\

while, unfortunately on the other hand

\begin{itemize}

\item The mappings in this natural classical context (8.2) are typically {\it not} surjective.
In other words, linear or nonlinear PDEs in (2.1) typically {\it cannot} be expected to have
{\it classical} solutions $U \in {\cal C}^m ( \Omega )$, for arbitrary continuous right hand
terms $f \in {\cal C}^0 ( \Omega )$. \\
Furthermore, it can often happen that nonclassical solutions have a major applicative interest,
thus they have to be sought out beyond the classical framework in (8.2).

\end{itemize}

This is, therefore, how we are led to the {\it necessity} to consider {\it generalized
solutions} $U$ for PDEs of type (2.1), that is, solutions $U \notin {\cal C}^m ( \Omega )$,
which therefore are no longer classical. This means that the natural classical mappings (8.2)
must in certain suitable ways be {\it extended} to {\it commutative diagrams} \\

\bigskip
\begin{math}
\setlength{\unitlength}{0.2cm}
\thicklines
\begin{picture}(60,21)

\put(10,19){${\cal C}^m ( \Omega )$}
\put(27,21){$T ( x, D )$}
\put(18,19.5){\vector(1,0){26.5}}
\put(47,19){${\cal C}^0 ( \Omega )$}
\put(0,11){$(8.3)$}
\put(12,17){\vector(0,-1){12}}
\put(13.5,11){$\subseteq$}
\put(49,17){\vector(0,-1){12}}
\put(45.7,11){$\subseteq$}
\put(10,2){$~~{\cal X}$}
\put(16,2.6){\vector(1,0){29.5}}
\put(47,2){$~~{\cal Y}$}
\put(29,-0.5){${\cal T}$}

\end{picture}
\end{math} \\

which are expected to have certain kind of {\it surjectivity} type properties, such as for
instance \\

(8.4)~~~ $ {\cal C}^0 ( \Omega ) ~\subseteq~ T ( {\cal X} ) $ \\

We conclude with a few comments :

\begin{itemize}

\item
Traditionally, ever since Hilbert and Sobolev, starting before WW II, {\it functional analysis}
has been used in solving PDEs, and suitable uniform topologies are defined on the domains and
ranges of the corresponding partial differential operators $T ( x, D )$. Thus these operators
obtain certain continuity properties. Then the extensions ${\cal X}$ and ${\cal Y}$ in (8.3)
are defined as the {\it completions} in these uniform topologies of the domains and ranges of
$T ( x, D )$, respectively. Finally, the continuity properties of $T ( x, D )$ may allow the
construction of suitable extensions ${\cal T}$ which would give the commutative diagrams (8.3),
and also satisfy some version of the surjectivity property (8.4), see for details
Oberguggenberger \& Rosinger [chap. 12, pp. 237-262].

\item
Since the 1960s, the {\it algebraic} method in 46F30 can alternatively be used especially in
the case of {\it nonlinear} partial differential operators $T ( x, D )$. In this respect, large
classes of {\it differential algebras of generalized functions} containing the Schwartz
distributions were constructed as the sought after extensions ${\cal X}$ and ${\cal Y}$ in
(8.3). \\
The most general classes of such algebras were introduced and used in Rosinger [1-17], starting
with the 1960s. Later, in the 1980s, a particular class of such algebras was introduced in
Colombeau, and it has known a certain popularity. However, due to the specific polynomially
limiting {\it growth conditions} required in the construction of Colombeau algebras, their use
in the study, for instance, of Lie group symmetries of PDEs, or singularities in General
Relativity is limited, since in both cases one may have to deal with transformation whose
growth can be arbitrary. In this way, such transformations cannot be accommodated within the
Colombeau algebras of generalized functions. On the other hand, arbitrary smooth
transformations and operations can easily be dealt with in the much larger classes of algebras
of generalized functions introduced earlier in Rosinger [1-17].

\item
The {\it order completion} method, introduced and developed in 1994 in Oberguggenberger \&
Rosinger, and further improved in Anguelov \& Rosinger, constructs the extensions ${\cal X}$
and ${\cal Y}$ in (8.3) as the {\it Dedekind order completion} of spaces naturally associated
with the partial differential operators $T ( x, D )$, and the spaces ${\cal C}^m ( \Omega )$
and ${\cal C}^0 ( \Omega )$ in (8.2).

\end{itemize}

Related to the advantages of the order completion method in solving nonlinear PDEs let us
mention here in short the following. \\

Neither the functional analytic, nor the algebraic methods can so far come anywhere near to
solve nonlinear PDEs of the generality of those in (2.1), let alone systems of such nonlinear
PDEs together with associated initial and/or boundary value problems. \\
In fact, the functional analytic methods are still subjected to the celebrated 1957 Hans Lewy
impossibility which they are nowhere near to manage to overcome, even in the general smooth
coefficient linear case. \\

As far as the algebraic method is concerned, it has among others come quite near to the
solution of the Lewy impossibility, see Colombeau, and for a short respective account
Rosinger [7]. \\
Further powerful results in solving various classes of nonlinear PDEs, not treated so far by
the functional analytic method, can be found in Rosinger [5-10], Colombeau , or
Oberguggenberger. \\
Among such results is the {\it global} solution of arbitrary analytic nonlinear systems of
PDEs, when considered with analytic non-characteristic Cauchy initial values, mentioned in
section 4 above. The respective generalized solutions obtained are {\it analytic} functions,
except for {\it closed} and {\it nowhere dense} subsets $\Gamma$ of the domains $\Omega$ of
definition of the given PDEs. In addition, these subsets $\Gamma$ can also be chosen to have
{\it zero} Lebesgue measure. \\

On the other hand, the order completion method introduced and developed in Oberguggenberger \&
Rosinger, and further improved in Anguelov \& Rosinger, as far as the {\it regularity} of
solutions is concerned, can not only deliver global solutions for systems of nonlinear PDEs of
the generality of those in (2.1), but it can also obtain a universal or blanket type
independent {\it regularity} result for such solution, namely, prove that they can be
assimilated with usual measurable functions, or even with {\it Hausdorff continuous}
functions. \\

Clearly, as one of the consequences of solving nonlinear systems of PDEs of the generality of
those in (2.1), the order completion method in Oberguggenberger \& Rosinger and in Anguelov \&
Rosinger is the only one so far which fully manages to overcome the Lewy impossibility.
Furthermore, it does so with a {\it large} nonlinear margin. \\ \\

{\bf 9. Order Completion Abolishes the Dichotomy "Linear \\
\hspace*{0.5cm} versus Nonlinear"} \\

The {\it dichotomy linear versus nonlinear} relating to equations or operators is in its
essence an issue of {\it algebra}, and more specifically, of {\it vector space} structures. In
this way, it is present {\it both} in the functional analytic {\it and} algebraic methods for
solving PDEs. And needless to say, dealing with the nonlinear case proves to be incomparably
more difficult than it is with the linear one. Consequently, and unfortunately, the presence
of this dichotomy is one of the {\it major disadvantages} of both the functional analytic and
algebraic methods in solving PDEs, even if by now it is taken so much for granted that no
attempt is made to abolish it. \\

On the other hand, order structures are of a more {\it basic} type than the algebraic ones. \\
Consequently, if instead of algebraic structures we consider order structures on the spaces of
smooth functions on which the partial differential operators act naturally when seen in the
classical context, see for instance (8.2), then these order structures - being more basic
than algebra - {\it can no longer distinguish} between the linearity or nonlinearity of such
partial differential operators. In this way, the traditional dichotomy between the linear and
nonlinear is simply set aside, and then the only problem left is whether indeed one can solve
PDEs in the completion of such order structures. \\
Fortunately, as shown in Oberguggenberger \& Rosinger, and also in Anguelov \& Rosinger, such
a solution is possible and useful. \\

Here of course, one may think that it may help if the respective partial differential
operators are monotonous. And then one may be concerned that all we managed to do was simply
to get rid of the dichotomy between linear and nonlinear, so that instead, now we have to face
the dichotomy monotonous versus arbitrary partial differential operators, with the latter
being quite likely again far more difficult to deal with, than the former. \\
Such a particular approach in which the dichotomy monotonous versus non-monotonous prevails
has recently been pursued in Carl \& Heikkil\"{a}, for instance. \\

This however is clearly {\it not} the way in the method in Oberguggenberger \& Rosinger and in
Anguelov \& Rosinger Indeed, although in this method order structures and order completions are
essentially used in solving PDEs, one need {\it not} require a priori any sort of monotonicity
property related to the equations solved. Certainly, the generality of nonlinear PDEs in
(2.1) or of systems of such nonlinear PDEs illustrates the fact that the respective equations
are {\it not} supposed to satisfy any a priori monotonicity conditions whatsoever. \\

What happens is very simple in fact, and is {\it similar} with the way the operators
$T ( x, D )$ and ${\cal T}$ in (8.3) acquire continuity type properties when the functional
analytic method is used in the construction of such commutative diagrams. Indeed, with the
functional analytic method, when one starts, say, with the natural mappings (8.2), the
uniform topologies one considers on the classical domains and ranges of the operators in order
to obtain commutative diagrams (8.3) are not arbitrary, but typically are related to the
respective operators, see section 6 above, or Oberguggenberger \& Rosinger [chap. 12]. \\

The same happens when order completion is used in Oberguggenberger \& Rosinger and in Anguelov
\& Rosinger for the construction of commutative diagrams (8.3). More precisely, the order
structures on the spaces of smooth functions on which the partial differential operators
naturally act, see for instance (8.2), will typically be defined {\it dependent} to a
certain extent on these operators. In this way, the respective operators, no matter how
arbitrary within the class of those in (2.1), will nevertheless {\it become} monotonous, thus
so much easier to deal with. \\

Such a procedure obviously cannot be imitated within algebra, since a nonlinear operator cannot
in general define a vector space structure in which it would become linear. \\
Furthermore, such a procedure obviously goes far beyond the approach in Carl \& Heikkil\"{a},
for instance, where one starts with given natural order relations on {\it both} the domains and
ranges of ODEs or PDEs, and then severely restricts oneself only to those rather small classes
of equations whose associated operators, or rather, inverse operators, are monotonous in the a
priori given orders. \\ \\

{\bf 10. The Hidden Power of Methods Based on Partial Orders} \\

As it happens, there is a rather widespread perception in mathematics that order structures
are too simple, and thus powerless, especially in analysis, therefore, they can deliver less
than algebra, which on its turn, can deliver less than topology. \\
Accordingly, since the emergence of functional analytic methods in the solution of PDEs at the
beginning of the 20th century, with the respective wealth of topologies on a variety of spaces
of functions, the perception prevails that there simply cannot be any other more powerful
methods in present day mathematics which could deal with such equations. \\

In view of such a perception it may appear surprising to see results such as in Anguelov \&
Rosinger, and the earlier ones in Oberguggenberger \& Rosinger, results obtained through order
structures, and which so far could not be approached anywhere near by functional analytic
methods. \\
Indeed, this method - based on order completion - does solve systems of PDEs of the nonlinear
generality of those in (2.1), together with associated initial and/or boundary value
problems, and furthermore, delivers for them global solutions which can be assimilated with
usual measurable or even Hausdorff continuous functions. It is in this way that the order
completion method is not only unprecedented, but it may also look rather strange in view of
the mentioned perception in mathematics related to order structures. \\

Therefore, one should indeed address the apparent {\it secret of the power} of the order
completion method in solving such large classes of systems of nonlinear PDEs, together with
their associated initial and/or boundary value problems. And one can do so by questioning the
mentioned perception. This can perhaps best be done by the presentation of certain classical -
even if less well known - examples which illustrate the power of order structures in yielding
what usually are called {\it deep theorems}. \\

In this regard a rather impressive, yet less well known fact is given by the 1936 "Spectral
Theorem" of Freudenthal, see Luxemburg \& Zaanen [chap. 6].  \\
Let us recall here in short some of its rather deep consequences. \\
The mentioned "Spectral Theorem" is a theorem about partially ordered structures, and it was
proved by Freudenthal exclusively in terms of such structures. Yet, what is of special
relevance is that by suitable particularizations, one can obtain from it the following three
results which are in fields as diverse as Operator Theory, Measure Theory and linear PDEs :

\begin{itemize}

\item the celebrated spectral representation theorem for normal operators in Hilbert spaces,

\item the highly nontrivial Radon-Nikodym theorem in measure theory, and

\item the Poisson formula for harmonic functions in an open circle.

\end{itemize}

{~} \\

{\bf 11. Conclusions} \\

The unprecedented power of the order completion method in solving very general systems of
nonlinear PDEs and the associated initial and/or boundary value problems stems from two facts :

\begin{itemize}

\item Partial orders are more basic mathematical structures than algebra or topology. And being
more basic than algebra, partial orders do {\it not} distinguish between linear and nonlinear
equations, operators, and so on. Consequently, partial orders treat the linear and nonlinear
cases in the same manner. Functional analytic method clearly cannot do the same.

\item When using order completion for solving PDEs, one need {\it not} assume any monotonicity
properties of the respective equations.

\end{itemize}

In the order completion method, the partial orders on the spaces of functions which are the
domains of definition of the partial differential operators considered are defined in relation
to these operators. This is similar to the way the topologies on such domains are defined, when
functional analytic methods are used. Further details in this regard can be found in
Oberguggenberger \& Rosinger [chap. 12, 13], or Anguelov \& Rosinger. \\

\end{document}